\newfont{\Bbb}{msbm10 scaled\magstephalf}
\documentclass[12pt,reqno]{amsart}
\usepackage{amsthm,amsfonts,amssymb,amsmath,amscd}
\usepackage{latexsym}
\usepackage{mathrsfs}
\usepackage{epsfig}
\usepackage{subfigure}
\usepackage{graphicx}
\usepackage{float}
\usepackage{dsfont}
\usepackage{amsthm}

\newtheorem{thm}{Theorem}[section]
\newtheorem{cor}[thm]{Corollary}
\newtheorem{Lemma}[thm]{Lemma}
\newtheorem{prop}[thm]{Proposition}
\newtheorem{ex}[thm]{Example}
\theoremstyle{definition}

\theoremstyle{remark}
\newtheorem{rem}[thm]{Remark}

\numberwithin{equation}{section}

\begin{document}
\title[Normal complex symmetric]{Normal complex symmetric weighted composition operators on the Hardy space}
\author{\bf Hang Zhou and Ze-Hua Zhou$^*$}
\address{\newline Hang Zhou\newline School of Mathematics, Tianjin University, Tianjin 300354, P.R. China.}
\email{cqszzs123@163.com}
\address{\newline Ze-Hua Zhou\newline School of Mathematics, Tianjin University, Tianjin 300354, P.R. China.}
\email{zehuazhoumath@aliyun.com;zhzhou@tju.edu.cn}
\date{}
\thanks{\noindent$^{*}$Corresponding author.
\\
The work was supported in part by the National Natural
Science Foundation of China (Grant Nos. 11771323; 11371276).}

\begin{abstract}
In this paper, we investigate the normal weighed composition operators $W_{\psi,\varphi}$ which is $\mathcal{J}-$symmetric, $\mathcal{C}_1-$symmetric and $\mathcal{C}_2-$symmetric on the Hardy space $H^2(\mathbb{D})$ respectively. Firstly, equivalent conditions of the normality of $\mathcal{C}_1-$symmetric and $\mathcal{C}_2-$symmetric weighted composition operators on $H^2(\mathbb{D})$ is given. Furthermore, the normal $\mathcal{J}-$symmetric, $\mathcal{C}_1-$symmetric and $\mathcal{C}_2-$symmetric weighted composition operators on $H^2(\mathbb{D})$ when $\varphi$ has an interior fixed point, $\varphi$ is of hyperbolic type or parabolic type are respectively investigated.
\end{abstract}
\maketitle
\renewcommand{\theequation}{\arabic{section}.\arabic{equation}}

\section{Introduction}
Let $\mathcal{B(H)}$ be the algebra of all bounded linear operators on a separable complex Hilbert space $\mathcal{H}$. A conjugation on $\mathcal{H}$ is an anti-linear operator $C:\mathcal{H}\rightarrow \mathcal{H}$ satisfying $\langle Cx,Cy \rangle=\langle y,x\rangle$ for all $x,y\in \mathcal{H}$ and $C^{2}=I$. An operator $T\in \mathcal{B(H)}$ is said to be complex symmetric if there exists a conjugation $C$ on $\mathcal{H}$ such that $T=CT^*C$. In this case, we say that $T$ is complex symmetric with the specific conjugation $C$. As it is already known, normal operator, Hankel operators, compressed Toeplitz operators and the Volterra operators are complex symmetric. Details on the complex symmetric operators were introduced by Garcia and  Putinar in \cite{GP1} and \cite{GP2} and the references therein.

Let $\mathbb{D}$ be the unit disk of the complex plane $\mathbb{C}$. Let $H(\mathbb{D})$ and $S(\mathbb{D})$ denote the collection of all analytic functions on $\mathbb{D}$ and all holomorphic self-maps of the unit disk respectively.  The $n$-th iterates of an analytic self-map $\varphi\in S(\mathbb{D})$, is denoted by $\varphi_n$ with $\varphi_0$ standing for the identity function, where $n=1,2,\cdots$. The one-to-one holomorphic function which maps $\mathbb{D}$ onto itself, called the \textit{M}$\ddot{o}$\textit{bius} transformation, denoted by $Aut(\mathbb{D})$, with the form $\lambda\varphi_a$, where $|\lambda|=1$ and $\varphi_a$ is the defined by
$\varphi_a(z)=\frac{a-z}{1-\overline{a}z}, a,z\in\mathbb{D}.$
The Hardy space, denoted by $H^2(\mathbb{D})$, is the collection of analytic functions $f(z)=\sum_{n=0}^\infty a_nz^n$ such that $$\|f\|^2=\sum_{n=0}^\infty|a_n|^2<\infty,z\in\mathbb{D},$$
alternatively, $$\|f\|^2=\sup_{0<r<1}\frac{1}{2\pi}\int_0^{2\pi}|f(re^{i\theta})|^2d\theta.$$
For each $w\in\mathbb{D}$, the reproducing kernel at $w$ is defined by $$K_w(z)=\frac{1}{1-\overline{w}z},z\in\mathbb{D}.$$ It is easy to check that for $f\in H^2(\mathbb{D}),\langle f,K_w \rangle=f(w).$ The normalized reproducing kernel at $w$ is defined by $k_w(z)=\frac{(1-|w|^2)^{\frac{1}{2}}}{1-\overline{w}z}.$

The composition operator $C_\varphi$ induced by $\varphi\in S(\mathbb{D})$ is defined as
$C_\varphi f(z)=f\circ\varphi(z),f\in H(\mathbb{D}),z\in\mathbb{D}.$
Given $\psi\in H(\mathbb{D})$, then the pointwise multiplication operator $M_\psi(f)=\psi\cdot f$ for all $f\in H(\mathbb{D})$ can be induced. Combining the composition operator $C_\varphi$ and the multiplication operator $M_\psi$, the weighted composition operator $W_{\psi,\varphi}$ is defined by $$W_{\psi,\varphi}f(z)=\psi(z)f(\varphi(z)),f\in H(\mathbb{D}).$$ We refer to the book \cite{CM} for more details about the composition operators on the spaces of analytic functions.

The study of complex symmetric (weighted) composition operators has a history of several years, see, e.g. \cite{BN1,GH,GW,GZ1,JKKL,N1,N2,NST}. Since normal operators are complex symmetric, it is natural to ask that are there normal complex symmetric (weighted) composition operators on $H^2(\mathbb{D})$? As it was proved in \cite{JKKL}, all $\mathcal{J}-$symmetric composition operators are normal (see, Corollary 3.10 in \cite{JKKL}), where $\mathcal{J}$ is the standard conjugation with $\mathcal{J}f(z)=\overline{f(\overline{z})}$ for $f\in H(\mathbb{D})$. Furthermore, Junga etc. in \cite{JKKL} presented the equivalent conditions for a $\mathcal{J}-$symmetric weighted composition operator to be normal (see, \cite{JKKL} Corollary 3.7).

In \cite{LK}, Lim and Khoi classified the weighted anti-linear conjugation $\mathcal{A}_{u,v}$ into two cases $\mathcal{C}_1$ and $\mathcal{C}_2$ and presented the equivalent conditions of the $\mathcal{C}_1-$symmetric and $\mathcal{C}_2-$symmetric $W_{\psi,\varphi}$ on $H_\gamma(\mathbb{D})$. In \cite{BN2}, Bourdon and Narayan investigate the normality of $W_{\psi,\varphi}$ on $H^2(\mathbb{D})$ when $\varphi$ has an interior fixed point or a Denjoy-Wolff point on $\overline{\mathbb{D}}$ respectively. Thus it is natural to find all normal $\mathcal{C}_1-$symmetric and $\mathcal{C}_2-$symmetric weighted composition operators.

Furthermore, arose by the open question raised by Noor in \cite{N2} that ``Does there exist a non-constant and non-automorphic symbol $\varphi$ for which $C_\varphi$ is complex symmetric but not normal on $H^2(\mathbb{D})$?'', we are also interested in another question: What are the normal complex symmetric $W_{\psi,\varphi}$ with $\varphi\in Aut(\mathbb{D})$ on $H^2(\mathbb{D})$?


In this paper, our work base on the known $\mathcal{J}-$symmetric, $\mathcal{C}_1-$\\symmetric and $\mathcal{C}_2-$symmetric weighted composition operators. Firstly, we present the equivalent conditions of normal $\mathcal{C}_1-$symmetric and $\mathcal{C}_2-$symmetric weighted composition operators on $H^2(\mathbb{D})$ according to the fixed point of the induced symbol $\varphi$. Furthermore, we investigate the normal $\mathcal{J}-$symmetric, $\mathcal{C}_1-$symmetric and $\mathcal{C}_2-$symmetric weighted composition operators on $H^2(\mathbb{D})$ when $\varphi$ has an interior fixed point, $\varphi$ is of hyperbolic type or parabolic type respectively.

It is noted that some results in this paper contains trivial but tedious calculations. We only show some essential details and omit the tedious part for the convenience of readers.

\section{Preliminaries}
\textbf{2.1 Cowen's formula for the adjoint of a linear-fractional composition operator}

For a non-constant linear-fractional self-map of the unit disk, Cowen in \cite{C2} established the important formula
$C_\varphi^*=M_g C_\sigma M_h^*,$
where the Cowen auxiliary functions are defined as follows:
\begin{align}\label{auxiliary functions}
g(z)=\frac{1}{-\overline{b}z+\overline{d}}, \sigma(z)=\frac{\overline{a}z+\overline{c}}{-\overline{b}z+\overline{d}},
h(z)=cz+d.\end{align}
\textbf{2.2 Normality of weighted composition operators}

In \cite{BN2} Bourdon and Narayan described all normal weighted composition operators $W_{\psi,\varphi}$ on $H^2(\mathbb{D})$, where $\varphi$ has an interior fixed point. They also presented an equivalent condition of the normality of $W_{\psi,\varphi}$ if $\varphi$ has a Denjoy-Wolff point on $\overline{\mathbb{D}}$. The results is of vital importance for our work and thus we show them as follows.
\begin{prop}\label{normal+interior fixed point}
\textup{[}\cite{BN2},Theorem 10\textup{]} Suppose that $\varphi$ has a fixed point $p\in\mathbb{D}$. Then $W_{\psi,\varphi}$ is normal if and only if $\psi=\psi(p)\frac{K_p}{K_p\circ\varphi}$ and $\varphi=\varphi_p\circ(\delta\varphi_p),$
where $\varphi_p(z)=\frac{p-z}{1-\overline{p}z}$ and $\delta\in\mathbb{C}$ with $|\delta|\leq1$.
\end{prop}

\begin{prop}\label{normal+D-W point on the boundary}
\textup{[}\cite{BN2},Theorem 12\textup{]} Suppose that $\varphi(z)=\frac{az+b}{cz+d}$ is a linear-fractional self-map of the unit disk and $\psi=K_{\sigma(0)},$ where $\sigma(z)=\frac{\overline{a}z-\overline{c}}{-\overline{b}z+\overline{d}}.$ Then $W_{\psi,\varphi}$ is normal if and only if
$$\frac{|d|^2}{|d|^2-|b|^2-(\overline{b}a-\overline{d}c)z}C_{\sigma\circ\varphi}=\frac{|d|^2}{|d|^2-|c|^2-(\overline{b}d-c\overline{a})z}C_{\varphi\circ\sigma}.$$
\end{prop}

\begin{rem}
It is easily checked that $$\frac{|d|^2}{|d|^2-|b|^2-(\overline{b}a-\overline{d}c)z}C_{\sigma\circ\varphi}=\frac{|d|^2}{|d|^2-|c|^2-(\overline{b}d-c\overline{a})z}C_{\varphi\circ\sigma}$$ is equivalent with
\begin{align}\label{normal+D-W point on the boundary+condition}
|\varphi(0)|=|\sigma(0)| \quad\textup{and}\quad \varphi\circ\sigma=\sigma\circ\varphi.
\end{align}(see, also, [\cite{L}, Proposition 4.6])
\end{rem}

Bourdon and Narayan also showed that if $\varphi$ is a linear-fractional self-map of parabolic type and $\psi=K_{\sigma(0)},$ then $W_{\psi,\varphi}$ is normal. Furthermore, they showed that no hyperbolic non-automorphic linear-fractional map can induce a normal weighted composition operator under the condition of Theorem \ref{normal+D-W point on the boundary} (see, [\cite{BN2}, Proposition 13] and the remarks below).

\textbf{2.3 Complex symmetric weighted composition operators on the Hardy space}

In \cite{LK} Lim and Khoi presented the equivalent conditions of \\$\mathcal{J}-$symmetric, $\mathcal{C}_1-$symmetric and $\mathcal{C}_2-$symmetric weighted composition operators $W_{\psi,\varphi}$ on $H_\gamma(\mathbb{D})$, which obviously contains the results on $H^2(\mathbb{D})$. Since our work is based on their results, we show them as follows. The first theorem can be obtained by [\cite{LK},Proposition 2.3].

\begin{prop}\label{J-symmetric}
$W_{\psi,\varphi}$ is $\mathcal{J}-$symmetric if and only if for each $z\in\mathbb{D}$, $\psi(z)=\frac{b}{1-a_o z}$ and $\varphi(z)=a_0+\frac{a_1z}{1-a_0z},$
where $a_0\in\mathbb{D},a_1\in\mathbb{D}$ and $b\in\mathbb{C}$.
\end{prop}

Let $u\in H(\mathbb{D})$ and $v\in S(\mathbb{D})$, for each $f\in H^2(\mathbb{D})$ and $z\in\mathbb{D}$, the weighted anti-linear operator $\mathcal{A}_{u,v}$ on $H^2(\mathbb{D})$ is defined by
$$\mathcal{A}_{u,v}=u(z)\overline{f(\overline{v(z)})},$$
which is a generalization of the standard conjugation $\mathcal{J}$. Observe that $W_{u,v}=\mathcal{A}_{u,v}\mathcal{J}$ and $\mathcal{A}_{u,v}=W_{u,v}\mathcal{J}.$ It is proved (see, Theorem 2.11 in \cite{LK}) that $\mathcal{A}_{u,v}:H^2(\mathbb{D})\rightarrow H^2(\mathbb{D})$ is a conjugation if and only if it has either of the following form:(i) there exists $\alpha,\lambda\in\partial\mathbb{D}$ such that for all $z\in\mathbb{D},$
\begin{align}\label{C_1}u(z)=\lambda \quad\textit{and}\quad v(z)=\alpha z.\end{align}
(ii) there exist $\alpha\in\mathbb{D}\setminus\{0\}$ and $\lambda\in\partial\mathbb{D}$ such that for all $z\in\mathbb{D}$,
\begin{align}\label{C_2}u(z)=\lambda k_\alpha(z) \quad\textit{and}\quad v(z)=\frac{\overline{\alpha}}{\alpha}\frac{\alpha-z}{1-\overline{\alpha}z}.\end{align}
As it is noted in \cite{LK}, $\mathcal{C}_1$ denotes $\mathcal{A}_{u,v}$ with $u$ and $v$ of the form \eqref{C_1} and $\mathcal{C}_2$ denotes $\mathcal{A}_{u,v}$ with $u$ and $v$ of the form \eqref{C_2}. Then the two theorems in the following have their roots in [\cite{LK},Theorem 3.1] and [\cite{LK},Theorem 3.2].

\begin{prop}
$W_{\psi,\varphi}:H^2(\mathbb{D})\rightarrow H^2(\mathbb{D})$ is $\mathcal{C}_1-$symmetric if and only if for all $z\in\mathbb{D},$
\begin{align}\label{C_1-symmetric}\psi(z)=\frac{d}{1-\alpha c_0z} \quad\textit{and}\quad \varphi(z)=c_0+\frac{c_1z}{1-\alpha c_0z},\end{align}
where $c_0,c_1\in\mathbb{D}$ and $d\in\mathbb{C}$.
\end{prop}

\begin{prop}
$W_{\psi,\varphi}:H^2(\mathbb{D})\rightarrow H^2(\mathbb{D})$ is $\mathcal{C}_2-$symmetric if and only if for all $z\in\mathbb{D},$
$$\psi(z)=\frac{d(c_0^2-\alpha c_1)}{c_0^2-\alpha c_1-(c_1-c_2)z}$$ and
\begin{align}\label{C_2-symmetric}
\varphi(z)=\frac{\alpha(\overline{\alpha}c_0^2-c_1)-(|\alpha|^2c_1-c_2)z}{\overline{\alpha}(c_0^2-\alpha c_1)-\overline{\alpha}(c_1-c_2)z},\end{align}
where $c_0,c_1,c_2,d\in\mathbb{C}.$
\end{prop}

\section{Complex symmetric $W_{\psi,\varphi}$ with $\varphi\in Aut(\mathbb{D})$}
In \cite{LK}, the authors presented the explicit form of the $\mathcal{J}-$symmetric weighted composition operator $W_{\psi,\varphi}$ with $\varphi\in Aut(\mathbb{D}).$

\begin{Lemma}\label{J-symmetric+aut}

Let $\varphi(z)=a_0+\frac{a_1z}{1-a_0z}$ be the form in Theorem \ref{J-symmetric}, where $a_0,a_1\in\mathbb{D}.$ Then $\varphi\in Aut(\mathbb{D})$ if and only if\\
(i) there exists $\gamma\in\mathbb{D}\setminus\{0\}$ such that $\varphi(z)=\frac{\overline{\gamma}}{\gamma}\frac{\gamma-z}{1-\overline{\gamma}z},$
where $\gamma=\frac{a_1+1}{a_0},z\in\mathbb{D}.$\\
(ii) there exist $\beta\in\partial\mathbb{D}$ such that $\varphi(z)=\beta z,$
where $\beta=a_1,z\in\mathbb{D}.$
\end{Lemma}
\textit{Proof.} The necessity is from \textup{[}\cite{LK}, Proposition 2.8\textup{]} and the sufficiency is obviously checked.

However, the explicit forms of the $\mathcal{C}_1-$symmetric and $\mathcal{C}_2-$symmetric $W_{\psi,\varphi}$ with $\varphi\in Aut(\mathbb{D})$ are still unknown. Since it is essential for our study in this paper, we prove them in a similar way with [\cite{LK}, Proposition 2.8].

\begin{Lemma}\label{C_1-symmetric+aut}
Let $\varphi$ be the form \eqref{C_1-symmetric}. Then $\varphi\in Aut(\mathbb{D})$ if and only if\\
(i) there exists $\gamma\in\mathbb{D}\setminus\{0\}$ such that
\begin{align}\label{C_1-symmetric+aut+form}\varphi(z)=\frac{\overline{\gamma}}{\gamma\alpha}\frac{\gamma-z}{1-\overline{\gamma}z},\end{align}
where $\gamma=\frac{c_0}{\alpha c_0^2-c_1},z\in\mathbb{D}.$\\
(ii) there exist $\beta\in\partial\mathbb{D}$ such that $\varphi(z)=\beta z,$
where $\beta=c_1,z\in\mathbb{D}.$
\end{Lemma}
\textit{Proof.} Since $\varphi\in Aut(\mathbb{D}),$ there exists $\gamma\in\mathbb{D}\setminus\{0\},\beta\in\partial\mathbb{D}$ such that \begin{align}\label{3.1}\varphi(z)=\beta\frac{\gamma-z}{1-\overline{\gamma}z}.\end{align}
Equating $\varphi(z)=\frac{c_0+(c_1-\alpha c_0^2)z}{1-\alpha c_0z}$ from \eqref{C_1-symmetric} and \eqref{3.1}, we have that
$$c_0+(c_1-\alpha c_0^2-\overline{\gamma}c_0)z-\overline{\gamma}(c_1-\alpha c_0^2)z^2=\beta\gamma-(\beta+\beta\alpha\gamma c_0)z+\alpha\beta c_0z^2.$$
Comparing the constants and the coefficients of $z$ and $z^2$ respectively,
\begin{align}
&c_0=\beta\gamma,\label{3.2}\\
&c_1-\alpha c_0^2-\overline{\gamma}c_0=-\beta(1+\alpha\gamma c_0),\label{3.3}\\
&-\overline{\gamma}(c_1-\alpha c_0^2)=\alpha\beta c_0.\label{3.4}
\end{align}
If $\gamma=0,$ then the second form of this lemma is proved.

If $\gamma\in\mathbb{D}\setminus\{0\}$, we substitute \eqref{3.2} into \eqref{3.4} to obtain that $c_1=\frac{(|\gamma|^2-1)\alpha\beta^2\gamma}{\overline{\gamma}}.$ Then substituting $c_0$ and $c_1$ again into \eqref{3.3}, we obtain that $\beta^2\alpha\gamma(|\gamma|^2-1)=\beta\overline{\gamma}(|\gamma|^2-1).$ Since $|\gamma|\neq1$ and $\beta\neq0$, we have that $\beta=\frac{\overline{\gamma}}{\alpha\gamma}.$
Then we have
$$\varphi(z)=\frac{\beta\gamma+(\frac{(|\gamma|^2-1)\alpha\beta^2\gamma}{\overline{\gamma}}-\alpha(\beta\gamma)^2)z}{1-\alpha \beta\gamma z}=\frac{\overline{\gamma}}{\gamma\alpha}\frac{\gamma-z}{1-\overline{\gamma}z},$$ where we omit the tedious calculation. Moreover, an easy calculation shows that $\gamma=\frac{c_0}{\alpha c_0^2-c_1}.$ The converse part is obviously checked. This completes the proof.\hfill$\Box$

\begin{Lemma}\label{C_2-symmetric+aut}
Let $\varphi$ be the form \eqref{C_2-symmetric}. Then $\varphi\in Aut(\mathbb{D})$ if and only if\\
(i) there exists $\gamma\in\mathbb{D}\setminus\{0\}$ such that
\begin{align}\label{C_2-symmetric+aut+form}\varphi(z)=\frac{|\alpha|^2-\alpha\overline{\gamma}}{\overline{\alpha}\gamma-|\alpha|^2}\frac{\gamma-z}{1-\overline{\gamma}z},z\in\mathbb{D},\end{align}
where $$\gamma=\frac{\alpha(\overline{\alpha}c_0^2-c_1)}{|\alpha|^2c_1-c_2}.$$\\
(ii) $\varphi(z)=z,z\in\mathbb{D}.$
\end{Lemma}
\textit{Proof.} Since $\varphi\in Aut(\mathbb{D}),$ there exists $\gamma\in\mathbb{D}\setminus\{0\},\beta\in\partial\mathbb{D}$ such that \begin{align}\label{3.5}\varphi(z)=\beta\frac{\gamma-z}{1-\overline{\gamma}z}.\end{align}
Equating $\varphi(z)=\frac{\alpha(\overline{\alpha}c_0^2-c_1)-(|\alpha|^2c_1-c_2)z}{\overline{\alpha}(c_0^2-\alpha c_1)-\overline{\alpha}(c_1-c_2)z}$ from \eqref{C_2-symmetric} and \eqref{3.5} and comparing the constants and the coefficients of $z$ and $z^2$ respectively, we have that
\begin{align}
&\overline{\gamma}(|\alpha|^2c_1-c_2)=\beta\overline{\alpha}(c_1-c_2),\label{3.6}\\
&\alpha(\overline{\alpha}c_0^2-c_1)=\beta\gamma\overline{\alpha}(c_0^2-\alpha c_1),\label{3.7}\\
&\overline{\gamma}\alpha(\overline{\alpha}c_0^2-c_1)+(|\alpha|^2c_1-c_2)=\overline{\alpha}\beta\gamma(c_1-c_2)+\beta\overline{\alpha}(c_0^2-\alpha c_1).\label{3.8}
\end{align}
If $\gamma=0,$ then we get $c_1=c_2$, $c_0^2=\frac{c_1}{\overline{\alpha}}$ since $\alpha,\beta\neq0$ and
\begin{align}\label{3.9}|\alpha|^2c_1-c_2=\beta\overline{\alpha}(c_0^2-\alpha c_1).\end{align}
Putting $c_1=c_2$, $c_0^2=\frac{c_1}{\overline{\alpha}}$ into \eqref{3.9}, we have that $\beta=-1.$
It follows that $\varphi(z)=\frac{\alpha(\overline{\alpha}\frac{c_1}{\overline{\alpha}}-c_1)-(|\alpha|^2-1)c_1z}{\overline{\alpha}(\frac{c_1}{\overline{\alpha}}-\alpha c_1)}=z,$ which implies the second part of this lemma.

If $\gamma\in\mathbb{D}\setminus\{0\}$, then combining \eqref{3.6}, \eqref{3.7} and \eqref{3.8}, we have that
$\frac{|\alpha|^2c_1-c_2}{\overline{\alpha}c_0^2-c_1}=\frac{\alpha}{\gamma}$ and
$\frac{|\alpha|^2c_1-c_2}{c_0^2-\alpha c_1}=\overline{\alpha}\beta,$
which implies that
\begin{align}
&c_0^2=\frac{1}{\overline{\alpha}}\frac{|\alpha|^2\beta\gamma-\alpha}{\beta\gamma-\alpha}c_1 \label{C_2-symmetric+aut+c_0}\\
&c_2=(1-\frac{\alpha\overline{\gamma}}{\overline{\alpha}}\frac{|\alpha|^2-1}{\beta\gamma-\alpha})c_1.
\label{C_2-symmetric+aut+c_2}
\end{align}
Putting $c_0^2$ and $c_1$ into $\frac{|\alpha|^2c_1-c_2}{c_0^2-\alpha c_1}=\overline{\alpha}\beta$, we get
\begin{align}
\beta=\frac{|\alpha|^2-\alpha\overline{\gamma}}{\overline{\alpha}\gamma-|\alpha|^2}.\label{C_2-symmetric+aut+beta}
\end{align}
Then we have
\begin{align*}
\varphi(z)&=\frac{\alpha(\overline{\alpha}c_0^2-c_1)-(|\alpha|^2c_1-c_2)z}{\overline{\alpha}(c_0^2-\alpha c_1)-\overline{\alpha}(c_1-c_2)z}\\
&=\frac{\alpha(\overline{\alpha}\frac{1}{\overline{\alpha}}\frac{|\alpha|^2\beta\gamma-\alpha}{\beta\gamma-\alpha}c_1-c_1)-(|\alpha|^2c_1-(1-\frac{\alpha\overline{\gamma}}{\overline{\alpha}}\frac{|\alpha|^2-1}{\beta\gamma-\alpha})c_1)z}{\overline{\alpha}(\frac{1}{\overline{\alpha}}\frac{|\alpha|^2\beta\gamma-\alpha}{\beta\gamma-\alpha}c_1-\alpha c_1)-\overline{\alpha}(c_1-(1-\frac{\alpha\overline{\gamma}}{\overline{\alpha}}\frac{|\alpha|^2-1}{\beta\gamma-\alpha})c_1)z}\\
&=\beta\frac{\gamma-z}{1-\overline{\gamma}z}.
\end{align*}
Moreover, an easy calculation shows that $\gamma=\frac{\alpha(\overline{\alpha}c_0^2-c_1)}{|\alpha|^2c_1-c_2}.$ The converse part is obviously checked. This completes the proof.\hfill$\Box$
\section{Normal $\mathcal{J}-$symmetric $W_{\psi,\varphi}$ with $\varphi\in Aut(\mathbb{D})$}
This first proposition in the following give an equivalent condition of $\mathcal{J}-$symmetric weighted composition operators to be normal, which is from Corollary 3.7 in \cite{JKKL}.

\begin{prop}\label{normal+J-symmetric}
$W_{\psi,\varphi}:H^2(\mathbb{D})\rightarrow H^2(\mathbb{D})$ is $\mathcal{J}-$symmetric and normal if and only if $\psi(z)=\frac{b}{1-a_o z}$, $\varphi(z)=a_0+\frac{a_1z}{1-a_0z}$ and
\begin{align}\label{normal+J-symmetric+condition}\Im a_0-|a_0|^2 \Im a_0+\Im \overline{a_0}a_1=0,\end{align}
where $a_0,a_1\in\mathbb{D},b\in\mathbb{C},z\in\mathbb{D}.$
\end{prop}

\begin{rem}
Note that the proposition above can also be proved by Proposition \ref{normal+D-W point on the boundary} with the equivalent condition \eqref{normal+D-W point on the boundary+condition}.
\end{rem}

In the following, we investigate the normal $\mathcal{J}-$symmetric weighted composition operators with automorphism symbols.
\begin{cor}
Suppose that $\varphi$ satisfies the hypothesis in Proposition \ref{normal+J-symmetric}. Then $\varphi\in Aut(\mathbb{D}) $ if and only if it has either of the following forms:\\
(i) there exists $\alpha\in\mathbb{D}\setminus\{0\}$ and $\beta\in\partial\mathbb{D}$ such that $\varphi(z)=\frac{\overline{\alpha}}{\alpha}\frac{\alpha-z}{1-\overline{\alpha}z},z\in\mathbb{D}.$\\
(ii) there exist $\beta\in\partial\mathbb{D}$ such that $\varphi(z)=\beta z,z\in\mathbb{D}.$
\end{cor}

\textit{Proof.} By Lemma \ref{J-symmetric+aut}, we are only supposed to check \eqref{normal+J-symmetric+condition}.
If there exists $\alpha\in\mathbb{D}\setminus\{0\}$ and $\beta\in\partial\mathbb{D}$ such that $\varphi(z)=\frac{\overline{\alpha}}{\alpha}\frac{\alpha-z}{1-\overline{\alpha}z},$
where $a_0=\alpha\beta$ and $a_1=\frac{\beta^2\alpha(|\alpha|^2-1)}{\overline{\alpha}}.$ Substituting $a_0$ and $a_1$ into \eqref{normal+J-symmetric+condition}, we have that
$$\Im\alpha\beta-|\alpha\beta|^2\Im\alpha\beta+\Im\overline{\alpha\beta}\frac{\alpha\beta^2(|\alpha|^2-1)}{\overline{\alpha}}=0.$$
Furthermore, the second part is trivial to be checked. This completes the proof.\hfill$\Box$

Observe that Proposition \ref{normal+D-W point on the boundary} doesn't require that $\varphi$ has an interior fixed point (In fact, some special $\varphi$ with interior fixed point and $\delta$ with $|\delta|\leq1$ are included in Proposition \ref{normal+D-W point on the boundary}, e.g. $p=0$, $\delta\neq1$ or $p\neq0$, $\delta=1$ or $p\neq0,\delta\neq1,|\delta|=1$).

In the following, we give three examples by Proposition \ref{normal+J-symmetric} when $\varphi$ has an interior fixed point, $\varphi$ is of hyperbolic type or parabolic type respectively.

\begin{ex}\label{normal+J-symmetric+interior fixed point}
Suppose that $\varphi$ has an interior fixed point $p\in\mathbb{D}$, $\varphi$ is nonconstant. Then $W_{\psi,\varphi}:H^2(\mathbb{D})\rightarrow H^2(\mathbb{D})$ is $\mathcal{J}-$symmetric and normal if and only if
$$\psi(z)=\gamma\frac{1-p^2}{1-p^2\delta+p(\delta-1)z},\varphi(z)=\frac{p(1-\delta)+(\delta-p^2)z}{1-p^2\delta+p(\delta-1)z},$$
where $\gamma=\psi(p)\in\mathbb{C}, p\in(-1,1),z\in\mathbb{D}.$
\end{ex}
\textit{Proof.} Suppose that $\varphi,\psi$ have the form in Proposition \ref{normal+interior fixed point}, i.e.,
\begin{align}
&\psi(z)=\gamma\frac{1-|p|^2}{1-|p|^2\delta+\overline{p}(\delta-1)z},\label{normal+interior fixed point+psi}\\
&\varphi(z)=\frac{p(1-\delta)+(\delta-|p|^2)z}{1-|p|^2\delta+\overline{p}(\delta-1)z}\label{normal+interior fixed point+varphi}
\end{align}

Since $W_{\psi,\varphi}$ is $\mathcal{J}-$symmetric, $\psi,\varphi$ should coordinate with the form in Proposition \ref{J-symmetric}, i.e.,
\begin{align}
&\psi(z)=\frac{\psi(0)}{1-a_0z},\label{normal+J-symmetric+psi}\\
&\varphi(z)=\frac{(a_1-a_0^2)z+a_0}{1-a_0z}.\label{normal+J-symmetric+varphi}
\end{align}

Equating \eqref{normal+interior fixed point+psi} and \eqref{normal+J-symmetric+psi}, then comparing the constants and the coefficients of $z$ and $z^2$ respectively, we have that
\begin{align}\label{4.1}
\psi(0)=\frac{\gamma(1-|p|^2)}{1-|p|^2\delta}, a_0=\frac{\overline{p}(1-\delta)}{1-|p|^2\delta}.\end{align}
In fact, substituting $\psi(0)$ and $a_0$ of the forms above into \eqref{normal+J-symmetric+psi}, we get \eqref{normal+interior fixed point+psi}.

In the sequence, equating \eqref{normal+interior fixed point+varphi} and \eqref{normal+J-symmetric+varphi}, then comparing the constants and the coefficients of $z$ and $z^2$ respectively, we have that
\begin{align}
&a_0=\frac{p(1-\delta)}{1-|p|^2\delta},\label{4.2}\\
&a_1=a_0^2+\frac{a_0(|p|^2-\delta)}{\overline{p}(\delta-1)},\label{4.3}\\
&(a_1-a_0^2)(1-|p|^2\delta)+\overline{p}(\delta-1)a_0=(\delta-|p|^2)-p(1-\delta)a_0.\label{4.4}
\end{align}
Comparing \eqref{4.1} and \eqref{4.2}, we get $p\in\mathbb{R}.$ Furthermore, substituting \eqref{4.2} and \eqref{4.3} into \eqref{4.4}, we get
$(p-\overline{p})\delta+(p-\overline{p})\delta|p|^4-(p-\overline{p})2\delta|p|^2=0,$ which is trivial.
In fact, substituting $a_0=\frac{p(1-\delta)}{1-p^2\delta}$, $a_1=\delta\frac{(p^2-1)^2}{(1-p^2\delta)^2}$ into \eqref{normal+J-symmetric+varphi}, we also get \eqref{normal+interior fixed point+varphi}. This completes the proof.\hfill$\Box$

\begin{cor}
Suppose that $\varphi$ satisfies the hypothesis in Example \ref{normal+J-symmetric+interior fixed point}, then $\varphi\in Aut(\mathbb{D})$ if and only if it has either of the following forms:\\
(i) $\varphi(z)=\beta\frac{\alpha-z}{1-\overline{\alpha}z},$
where $\alpha=\frac{p(1-\overline{\delta})}{1-p^2\overline{\delta}}\in\mathbb{D},\beta=\frac{p^2-\delta}{1-p^2\delta}\in\partial\mathbb{D}$ and $\delta\in\partial\mathbb{D},z\in\mathbb{D}.$\\
(ii) $\varphi(z)=\delta z,$ where $\delta\in\partial\mathbb{D},z\in\mathbb{D}.$
\end{cor}
\textit{Proof.} By Lemma \ref{J-symmetric+aut}, if there exist $\alpha\in\mathbb{D}\setminus\{0\}$ and $\beta\in\partial\mathbb{D}$ such that $\varphi(z)=\frac{\overline{\alpha}}{\alpha}\frac{\alpha-z}{1-\overline{a}z},$
where $a_0=\alpha\beta$ and $a_1=\frac{\beta^2\alpha(|\alpha|^2-1)}{\overline{\alpha}}.$
By the proof of Example \ref{normal+J-symmetric+interior fixed point}, equating $a_0=\alpha\beta$ and $a_0=\frac{p(1-\delta)}{1-p^2\delta}$, we obtain that $\alpha=\frac{p(1-\delta)}{1-p^2\overline{\delta}}$ since $\beta=\frac{\overline{\alpha}}{\alpha}$. Likewise, equating the two forms of $a_1$, and substituting $\alpha=\frac{p(1-\delta)}{1-p^2\overline{\delta}}$ into it, we obtain that $\frac{1-\overline{\delta}}{1-p^2\overline{\delta}}=\frac{1-\delta}{p^2-\delta}.$
Hence, $\beta=\frac{\overline{\alpha}}{\alpha}=\frac{p^2-\delta}{1-p^2\delta}.$ Since $\beta\in\partial\mathbb{D},$ we get $\delta\in\partial\mathbb{D}.$ In fact, substituting $\alpha,\beta$ into $\varphi(z)=\frac{\overline{\alpha}}{\alpha}\frac{\alpha-z}{1-\overline{\alpha}z},$ we get \eqref{normal+interior fixed point+varphi} by an easy calculation.
Otherwise, if there exist $\beta\in\partial\mathbb{D}$ such that $\varphi(z)=\beta z,$ then obviously $\beta=\delta.$
This completes the proof.\hfill$\Box$

\begin{ex}\label{normal+J-symmetric+hyperbolic+aut}
There is no $\mathcal{J}-$symmetric and normal weighted composition operator $W_{\psi,\varphi}:H^2(\mathbb{D})\rightarrow H^2(\mathbb{D})$ if $\varphi$ is a hyperbolic automorphism linear-fractional self-map with Denjoy-Wolff point $1\in\partial\mathbb{D}$.
\end{ex}
\textit{Proof.} The result can be similarly proved by Example \ref{normal+C_1-symmetric+hyperbolic+aut}.

\begin{ex}\label{normal+J-symmetric+hyperbolic+non-aut}
There is no $\mathcal{J}-$symmetric and normal weighted composition operator $W_{\psi,\varphi}:H^2(\mathbb{D})\rightarrow H^2(\mathbb{D})$ if $\varphi$ is a hyperbolic non-automorphism linear-fractional self-map.
\end{ex}
\textit{Proof.} The result can be directly obtained by the remark below Proposition 13 in \cite{BN2}.

\begin{ex}\label{normal+J-symmetric+parabolic}
The weighted composition operator $W_{\psi,\varphi}:H^2(\mathbb{D})\rightarrow H^2(\mathbb{D})$ is $\mathcal{J}-$symmetric and normal when $\varphi$ is a parabolic linear-fractional self-map if and only if $\psi(z)=\frac{d}{1-a_0z}$ and $\varphi$ has either of the following forms:\\
(i) $\varphi(z)=\frac{(1-2a_0)z+a_0}{1-a_0z}, \Im a_0=|a_0|^2, a_0, z\in\mathbb{D}$.\\
(ii) $\varphi(z)=\frac{(1+2a_0)z+a_0}{1-a_0z}, \Im a_0=-|a_0|^2, a_0, z\in\mathbb{D}.$
\end{ex}
\textit{Proof.} Since $\varphi$ is a hyperbolic automorphism linear-fractional self-map, then we have either of the following assertions:

(i) $\varphi$ has the Denjoy-Wolff point $\zeta_1=1$ and $a_1=a_0-1.$ Hence, $\varphi(z)=\frac{(1-2a_0)z+a_0}{1-a_0z}$ and \eqref{normal+J-symmetric+condition} turns to be $\Im a_0=|a_0|^2.$

(ii) $\varphi$ has the Denjoy-Wolff point $\zeta_2=-1$ and $a_1=a_0+1.$ Hence, $\varphi(z)=\frac{(1+2a_0)z+a_0}{1-a_0z}$ and
\eqref{normal+J-symmetric+condition} turns to be $\Im a_0=-|a_0|^2.$ This completes the proof.\hfill$\Box$

\section{Normal $\mathcal{C}_1-$symmetric $W_{\psi,\varphi}$}
Recall that the anti-linear operator $\mathcal{A}_{u,v}$ on $H^2(\mathbb{D})$ is denoted by $\mathcal{C}_1$ if $u$ and $v$ have the forms \eqref{C_1}.

\begin{thm}\label{normal+C_1-symmetric}
Suppose that $\varphi$ is nonconstant, then $W_{\psi,\varphi}:H^2(\mathbb{D})\rightarrow H^2(\mathbb{D})$ is $\mathcal{C}_1-$symmetric and normal if and only if
$$\psi(z)=\frac{1}{1-\alpha c_o z}, \varphi(z)=\frac{(c_1-\alpha c_0^2)z+c_0}{1-\alpha c_0z}$$ and
\begin{align}\label{normal+C_1-symmetric+condition}(\overline{c_0}-\alpha c_0)(1-|c_0|^2)+\alpha c_0\overline{c_1}-\overline{c_0}c_1=0,\end{align}
where $c_0,c_1\in\mathbb{D},\alpha\in\partial\mathbb{D},z\in\mathbb{D}.$
\end{thm}
\textit{Proof.}
Firstly, since $W_{\psi,\varphi}$ is $\mathcal{C}_1-$symmetric, by \eqref{C_1-symmetric},
$$\psi(z)=\frac{d}{1-\alpha c_0z} \quad\textit{and}\quad \varphi(z)=\frac{(c_1-\alpha c_0^2)z+c_0}{1-\alpha c_0 z}.$$
Assume that $\varphi(z)=\frac{az+b}{cz+d}$ and $\psi(z)=K_{\sigma(0)}(z)=\frac{d}{cz+d}.$
Equating the two forms of $\varphi$, an easy calculation shows that $a=c_1-\alpha c_0^2,b=c_0$ and $c=-\alpha c_0.$
Similarly, equating the two forms of $\psi$, we get $d=\psi(0)=1.$
Furthermore, since $\sigma(z)=\frac{\overline{c_1-\alpha c_0^2}z+\overline{\alpha c_0}}{1-\overline{c_0}z}$,
\begin{align}\label{5.1}|\varphi(0)|=|c_0|=|\sigma(0)|.\end{align}
Another tedious calculation shows that for each $z\in\mathbb{D}$,
$$\varphi\circ\sigma(z)=\frac{(|c_1-\alpha c_0^2|-|c_0|^2)z+(c_0-|c_0|^2 c_0+\overline{\alpha c_0}c_1)}{(1-\alpha c_0z)(1-\overline{c_0}z)}$$
and
$$\sigma\circ\varphi(z)=\frac{(|c_1-\alpha c_0^2|-|c_0|^2)z+(\overline{\alpha c_0}-|c_0|^2\overline{\alpha c_0}+c_0\overline{c_1})}{(1-\alpha c_0z)(1-\overline{c_0}z)}.$$
It follows that $\sigma\circ\varphi=\varphi\circ\sigma$ if and only if
\begin{align}\label{5.2}(\overline{c_0}-\alpha c_0)(1-|c_0|^2)+\alpha c_0\overline{c_1}-\overline{c_0}c_1=0.\end{align}
By \eqref{normal+D-W point on the boundary+condition}, the normality follows. This completes the proof.\hfill$\Box$

\begin{cor}
Suppose that $\varphi$ satisfies the hypothesis in Theorem \ref{normal+C_1-symmetric}, then $\varphi\in Aut(\mathbb{D})$ if and only if it has either of the following forms:\\
(i) there exists $\gamma\in\mathbb{D}\setminus\{0\}$ and $\beta\in\partial\mathbb{D}$ such that
$\varphi(z)=\frac{\overline{\gamma}}{\gamma\alpha}\frac{\gamma-z}{1-\overline{\gamma}z}.$\\
(ii) there exist $\beta\in\partial\mathbb{D}$ such that $\varphi(z)=\beta z.$
\end{cor}
\textit{Proof.} By Lemma \ref{C_1-symmetric+aut}, we are only supposed to check \eqref{normal+C_1-symmetric+condition}, which obviously holds by an easy calculation. This completes the proof.\hfill$\Box$

In the following, we give three examples by Theorem \ref{normal+C_1-symmetric} when $\varphi$ has an interior fixed point, $\varphi$ is of hyperbolic type or parabolic type respectively.

\begin{ex}\label{normal+C_1-symmetric+interior fixed point}
Suppose that $\varphi$ has an interior fixed point $p\in\mathbb{D}$, $\varphi$ is nonconstant and $W_{\psi,\varphi}:H^2(\mathbb{D})\rightarrow H^2(\mathbb{D})$ is $\mathcal{C}_1-$symmetric. Then $W_{\psi,\varphi}$ is normal if and only if it has either of the following forms:\\
(i) $\psi(z)=\frac{\gamma(1-|p|^2)}{1-\alpha pz}$, $\varphi(z)=\frac{p-|p|^2z}{1-\alpha pz},$ where $\gamma=\psi(p)\in\mathbb{C},z\in\mathbb{D}.$\\
(ii) $\psi(z)=\frac{\gamma(1-\alpha p^2)}{1-\alpha p^2\delta+\alpha p(\delta-1)z},$
$\varphi(z)=\frac{p(1-\delta)+(\delta-\alpha p^2)z}{1-\alpha p^2\delta+\alpha p(\delta-1)z},$
where $\gamma=\psi(p)\in\mathbb{D},z\in\mathbb{D}.$
\end{ex}
\textit{Proof.} As what we do in Example \ref{normal+J-symmetric+interior fixed point}, equating the two forms of $\psi$, we get $\psi(0)=\frac{\gamma(1-|p|^2)}{1-|p|^2\delta}, c_0=\frac{\overline{p}(1-\delta)}{\alpha(1-|p|^2\delta)}.$ Further equating the two forms of $\varphi$, we get
\begin{align}
&c_0=\frac{p(1-\delta)}{(1-|p|^2\delta)},\label{5.3}\\
&c_1=\alpha c_0^2+\frac{\alpha c_0(|p|^2-\delta)}{\overline{p}(\delta-1)}\label{5.4},\\
&\delta-|p|^2-p(1-\delta)\alpha c_0=(1-|p|^2\delta)(c_1-\alpha c_0^2)+c_0\overline{p}(\delta-1)\label{5.5}.
\end{align}
Substituting \eqref{5.3} and \eqref{5.4} into \eqref{5.5}, by some tedious but trivial calculation, we have that
$$\delta(\alpha p-\overline{p})(|p|^4-2|p^2|+1)=0,$$ which only holds when $\overline{p}=\alpha p$ or $\delta=0$ since $p\in\mathbb{D}$.
If $\delta=0$, then $\psi(z)=\frac{\gamma(1-|p|^2)}{1-\alpha pz}$, $\varphi(z)=\frac{p-|p|^2z}{1-\alpha pz}.$ Observe that under this circumstance we must have $\alpha p=\overline{p}$ to equate
$\varphi(z)=\frac{p-|p|^2z}{1-\alpha pz}$ and $\varphi(z)=\frac{(c_1-\alpha c_0^2)z+c_0}{1-\alpha c_0z}.$
Therefore, $c_0=p,c_1=0$ in \eqref{C_1-symmetric}.
If $\alpha p=\overline{p},\delta\neq0$, then $c_0=\frac{p(1-\delta)}{1-\alpha p^2\delta}$,
$c_1=\alpha c_0^2+\frac{c_0(\alpha p^2-\delta)}{p(\delta-1)}$ in \eqref{C_1-symmetric}. Moreover, it is easily checked that $c_0\in\mathbb{D}$ and $c_1\in\mathbb{D}$ if and only if $\delta\neq1$. This completes the proof.\hfill$\Box$

\begin{cor}
Suppose that $\varphi$ satisfies the hypothesis in Example \ref{normal+C_1-symmetric+interior fixed point}. If $\varphi\in Aut(\mathbb{D}),$ then  $\varphi(z)=\frac{2p-(1+\alpha p^2)z}{1+\alpha p^2-2\alpha pz},z\in\mathbb{D}.$
\end{cor}
\textit{Proof.} Firstly consider $\varphi(z)=\frac{p-|p|^2z}{1-\alpha pz}$ with $c_0=p,c_1=0$ in \eqref{C_1-symmetric}. By Lemma \ref{C_1-symmetric+aut}, if $\varphi$ has the form \eqref{C_1-symmetric+aut+form}, then $p=\beta\gamma,\frac{(|\gamma|^2-1)p\alpha\beta}{\overline{\gamma}}=0,$ which implies that $p=0$ since $\alpha,\beta\neq0,|\gamma|\neq1,$ which further implies that $\varphi=0.$

Moreover, we consider $\varphi(z)=\frac{p(1-\delta)+(\delta-\alpha p^2)z}{1-\alpha p^2\delta+\alpha p(\delta-1)z}$ with $c_0=\frac{p(1-\delta)}{1-\alpha p^2\delta}$,
$c_1=\alpha c_0^2+\frac{c_0(\alpha p^2-\delta)}{p(\delta-1)}$ in \eqref{C_1-symmetric}.
Since $\beta=\frac{\overline{\gamma}}{\alpha\gamma}$ and $c_0=\beta\gamma,$ we have that $\gamma=\frac{p(1-\overline{\delta})}{1-|p|^2\overline{\delta}}$ (We still write $|p|^2$ here instead of $\alpha p^2$).
Equating $\varphi(z)=\frac{p(1-\delta)+(\delta-|p|^2)z}{1-|p|^2\delta+\overline{p}(\delta-1)z}$ and \eqref{C_1-symmetric+aut+form}, then comparing the constants and the coefficients of $z$ and $z^2$ respectively, we have that
\begin{align}
&\beta=\frac{|p|^2-\delta}{1-|p|^2\delta}\nonumber\\
&(|p|^2-\delta)(1-\overline{\delta})=(1-\delta)(1-|p|^2\delta),\label{5.6}\\
&|p|^2|1-\delta|^2(1-\delta)=|p|^2|1-\delta|^2(1-\overline{\delta}).\label{5.7}
\end{align}
Observe that \eqref{5.7} implies $\delta\in\mathbb{R}.$ Hence, by \eqref{5.6}, we have that $\delta=-1$ since $p\in\mathbb{D},$ which also implies that $\gamma=\frac{2p}{1+|p|^2}.$ Hence, $\varphi(z)=\frac{2p-(1+|p|^2)z}{1+|p|^2-2\overline{p}z}.$ This completes the proof.\hfill$\Box$

\begin{ex}\label{normal+C_1-symmetric+hyperbolic+aut}
There is no $\mathcal{C}_1-$symmetric and normal weighted composition operator $W_{\psi,\varphi}:H^2(\mathbb{D})\rightarrow H^2(\mathbb{D})$ if $\varphi$ is a hyperbolic automorphism linear-fractional self-map with Denjoy-Wolff point $1\in\partial\mathbb{D}$.
\end{ex}
\textit{Proof.} Since $\varphi$ is a hyperbolic automorphism linear-fractional self-map with Denjoy-Wolff point $1\in\partial\mathbb{D}$, we can assume that $$\varphi(z)=\frac{(r+1-t)z+r+t-1}{(r-t-1)z+r+t+1},$$ where $r=\frac{1}{\varphi'(1)}>1,t\in\mathbb{C}.$ Equating $\varphi(z)=\frac{(r+1-t)z+r+t-1}{(r-t-1)z+r+t+1}$ and \eqref{C_1-symmetric+aut+form}, then comparing the constants and the coefficients of $z$ and $z^2$ respectively, we have that
$$\gamma=\frac{r-t-1}{\alpha(r-t+1)}=\overline{\frac{\alpha(r+t-1)}{r+t+1}},$$ which implies that $r+t=1.$ Hence, $\gamma=0,$ which is impossible. This completes the proof.\hfill$\Box$

\begin{ex}\label{normal+C_1-symmetric+hyperbolic+non-aut}
There is no $\mathcal{C}_1-$symmetric and normal weighted composition operator $W_{\psi,\varphi}:H^2(\mathbb{D})\rightarrow H^2(\mathbb{D})$ if $\varphi$ is a hyperbolic non-automorphism linear-fractional self-map.
\end{ex}
\textit{Proof.} The result can be directly obtained by the remark below Proposition 13 in \cite{BN2}.

\begin{ex}\label{normal+C_1-symmetric+parabolic}
The weighted composition operator $W_{\psi,\varphi}:H^2(\mathbb{D})\rightarrow H^2(\mathbb{D})$ is $\mathcal{C}_1-$symmetric and normal when $\varphi$ is a parabolic linear-fractional self-map with the Denjoy-Wolff point $\zeta\in\partial\mathbb{D}$ if and only if
$\psi(z)=\frac{\zeta^2}{\zeta^2-c_0z}$, $\varphi(z)=\frac{(\zeta^2c_1-c_0^2)z+\zeta^2c_0}{\zeta^2-c_0z}$ and
$$(\zeta^2\overline{c_0}-c_0)(1-|c_0|^2)+c_0\overline{c_1}-\zeta^2\overline{c_0}c_1=0,$$
where $c_0,c_1\in\mathbb{D}, z\in\mathbb{D}.$
\end{ex}
\textit{Proof.} Since $\varphi$ is a hyperbolic automorphism linear-fractional self-map, then $\varphi$ has the Denjoy-Wolff point $\zeta=\frac{1+\alpha c_0^2-c_1}{2\alpha c_0}$ and $(c_1-\alpha c_0^2-1)^2=4\alpha c_0^2,$ where the latter one implies that $\alpha=\frac{1}{\zeta^2}$ by substituting the former one into it. Therefore, the result can be easily obtained by what we have observed above. This completes the proof.\hfill$\Box$

\section{Normal $\mathcal{C}_2-$symmetric $W_{\psi,\varphi}$}
Recall that the anti-linear operator $\mathcal{A}_{u,v}$ on $H^2(\mathbb{D})$ is denoted by $\mathcal{C}_2$ if $u$ and $v$ have the forms \eqref{C_2}.

Just as we calculate in Section 5, we should firstly equate the two forms of $\psi$ and $\varphi$.
Suppose that $\varphi(z)=\frac{az+b}{cz+d}$ and $\psi(z)=\frac{d}{cz+d}.$

Equating $\psi(z)=\frac{d}{cz+d}$ and $\psi(z)=\frac{\psi(0)(c_0^2-\alpha c_1)}{c_0^2-\alpha c_1-(c_1-c_2)z}$ in \eqref{C_2-symmetric}, then comparing the constants and the coefficients of $z$ and $z^2$ respectively, we have that
$\psi(0)=1$ and $c=\frac{d(c_2-c_1)}{c_0^2-\alpha c_1}.$
Therefore, $$\psi(z)=\frac{c_0^2-\alpha c_1}{c_0^2-\alpha c_1-(c_1-c_2)z}=\frac{d}{cz+d}.$$
Further equating $\varphi(z)=\frac{az+b}{cz+d}$ and $\varphi(z)=\frac{\alpha(\overline{\alpha}c_0^2-c_1)-(|\alpha|^2c_1-c_2)z}{\overline{\alpha}(c_0^2-\alpha c_1)-\overline{\alpha}(c_1-c_2)z}$ in \eqref{C_2-symmetric}, then comparing the constants and the coefficients of $z$ and $z^2$ respectively, we have that
\begin{align*}
a=\frac{-d(|\alpha|^2c_1-c_2)}{\overline{\alpha}(c_0^2-\alpha c_1)},b=\frac{d(|\alpha|^2c_0^2-\alpha c_1)}{\overline{\alpha}(c_0^2-\alpha c_1)},
c=\frac{d(c_2-c_1)}{c_0^2-\alpha c_1}.
\end{align*}
Substituting the expressions of $a,b,c$ above into $\varphi(z)=\frac{az+b}{cz+d}$, we obtain that
$$\varphi(z)=\frac{az+b}{cz+d}=\frac{\alpha(\overline{\alpha}c_0^2-c_1)-(|\alpha|^2c_1-c_2)z}{\overline{\alpha}(c_0^2-\alpha c_1)-\overline{\alpha}(c_1-c_2)z}.$$

\begin{thm}\label{normal+C_2-symmetric}
Suppose that $\varphi$ is nonconstant, then $W_{\psi,\varphi}:H^2(\mathbb{D})\rightarrow H^2(\mathbb{D})$ is $\mathcal{C}_2-$symmetric and normal if and only if
$$\psi(z)=\frac{c_0^2-\alpha c_1}{c_0^2-\alpha c_1-(c_1-c_2)z}$$ and
$$\varphi(z)=\frac{\alpha(\overline{\alpha}c_0^2-c_1)-(|\alpha|^2c_1-c_2)z}{\overline{\alpha}(c_0^2-\alpha c_1)-\overline{\alpha}(c_1-c_2)z}$$
and either of the followings holds:\\
(i) $|c_1-c_2|=|\overline{\alpha}c_0^2-c_1|=|c_0^2-\alpha c_1|\neq\frac{||\alpha|^2c_1-c_2|}{|\alpha|}$.\\
(ii) $|c_1-c_2|=|\overline{\alpha}c_0^2-c_1|=|c_0^2-\alpha c_1|=\frac{||\alpha|^2c_1-c_2|}{|\alpha|}$ and
$$\Im(\overline{A}-\overline{C})(\widetilde{A}+\widetilde{C})=0.$$
where $c_0,c_1,c_2,d\in\mathbb{C},\alpha\in\mathbb{D}\setminus\{0\},z\in\mathbb{D},$
$$A=(|\alpha|^2c_0^2-\alpha c_1)(\alpha\overline{c_0}^2-|\alpha|^2\overline{c_1}),C=\alpha(\overline{c_1}-\overline{c_2})(|\alpha|^2c_1-c_2),$$
$$\widetilde{A}=-\alpha(|\alpha|^2\overline{c_1}-\overline{c_2})(\overline{\alpha}c_0^2-c_1),
\widetilde{C}=|\alpha|^2(c_0^2-\alpha c_1)(\overline{c_1}-\overline{c_2}).$$
\end{thm}
\textit{Proof.}
Since $\sigma(z)=\frac{-(|\alpha|^2\overline{c_1}-\overline{c_2})z+\alpha(\overline{c_1}-\overline{c_2})}{-\overline{\alpha}(\alpha\overline{c_0}^2-\overline{c_1})z+\alpha(\overline{c_0}^2-\overline{\alpha c_1})},$ we have that
$|\sigma(0)|=\frac{|c_1-c_2|}{|c_0^2-\alpha c_1|},|\varphi(0)|=\frac{|\overline{\alpha}c_0^2-c_1|}{|c_0^2-\alpha c_1|}.$
Hence, $|\sigma(0)|=|\varphi(0)|$ if and only if $$|c_1-c_2|=|\overline{\alpha}c_0^2-c_1|.$$
Furthermore, some tedious calculations show that for each $z\in\mathbb{D}$,
$$\varphi\circ\sigma(z)=\frac{(D-B)z+(A-C)}{(\overline{C}-\overline{A})z+(E-B)},$$
and
$$\sigma\circ\varphi(z)=\frac{(D-E)z+(\widetilde{A}+\widetilde{C})}{(-\overline{\widetilde{A}}-\overline{\widetilde{C}})z+(E-B)},$$
where $$A=(|\alpha|^2c_0^2-\alpha c_1)(\alpha\overline{c_0}^2-|\alpha|^2\overline{c_1}),B=|\alpha|^2|\overline{\alpha}c_0^2-c_1|^2,$$
$$C=\alpha(\overline{c_1}-\overline{c_2})(|\alpha|^2c_1-c_2),D=||\alpha|^2c_1-c_2|^2,E=|\alpha|^2|c_0^2-\alpha c_1|^2,$$
$$\widetilde{A}=-\alpha(|\alpha|^2\overline{c_1}-\overline{c_2})(\overline{\alpha}c_0^2-c_1),
\widetilde{C}=|\alpha|^2(c_0^2-\alpha c_1)(\overline{c_1}-\overline{c_2}).$$
It follows that $\sigma\circ\varphi=\varphi\circ\sigma$ if and only if
\begin{align*}
&(D-B)(\overline{\widetilde{A}}+\overline{\widetilde{C}})=(D-E)(\overline{A}-\overline{C}),\\
&(A-C)(E-B)=(\widetilde{A}+\widetilde{C})(E-B),\\
&(D-B)(E-B)-(A-C)(\overline{\widetilde{A}}+\overline{\widetilde{C}})=\nonumber\\
&(D-E)(E-B)+(\tilde{A}+\tilde{C})(\overline{C}-\overline{A}).
\end{align*}
We assert that $B=E$ by a basic deduction. Then the conditions above turn to be
\begin{align}
&(D-B)(\overline{\widetilde{A}}+\overline{\widetilde{C}})=(D-B)(\overline{A}-\overline{C}),\label{6.1}\\
&\Im(\overline{A}-\overline{C})(\widetilde{A}+\widetilde{C})=0.\label{6.2}
\end{align}
If $D\neq B,$ by \eqref{6.1}, $A-C=\widetilde{A}+\widetilde{C},$ which implies that \eqref{6.2} always holds. If $D=B,$ then \eqref{6.1} always hold. Combining what we have observed above and \eqref{normal+D-W point on the boundary+condition}, the normality follows. This completes the proof.\hfill$\Box$

\begin{cor}\label{normal+C_2-symmetric+aut}
There is no $\varphi\in Aut(\mathbb{D})$ satisfing the hypothesis in Theorem \ref{normal+C_2-symmetric}.
\end{cor}
\textit{Proof.} Substituting $c_0^2=\frac{1}{\overline{\alpha}}\frac{|\alpha|^2\beta\gamma-\alpha}{\beta\gamma-\alpha}c_1$ and
$c_2=(1-\frac{\alpha\overline{\gamma}}{\overline{\alpha}}\frac{|\alpha|^2-1}{\beta\gamma-\alpha})c_1$ into
$|\overline{\alpha}c_0^2-c_1|=|c_0^2-\alpha c_1|,$ we have that
$|\alpha|(|\gamma|-1)(1-|\alpha|^2)=0,$ which implies that $|\gamma|=1,$ which is impossible. This completes the proof.\hfill$\Box$

In the following, we give three examples by Theorem \ref{normal+C_2-symmetric} when $\varphi$ has an interior fixed point, $\varphi$ is of hyperbolic type or parabolic type respectively.

\begin{ex}\label{normal+C_2-symmetric+interior fixed point}
Suppose that $\varphi$ has an interior fixed point $p\in\mathbb{D}$, $\varphi$ is nonconstant, then $W_{\psi,\varphi}:H^2(\mathbb{D})\rightarrow H^2(\mathbb{D})$ is $\mathcal{C}_2-$symmetric and normal if and only if
$\psi(z)=\frac{\psi(p)(1-p^2)}{1-p^2\delta+p(\delta-1)z},$
$\varphi(z)=\frac{p(1-\delta)+(\delta-p^2)z}{1-p^2\delta+p(\delta-1)z},$
where $p\in(-1,1),z\in\mathbb{D}.$
\end{ex}
\textit{Proof.}
To investigate the situation of $\varphi$ with an interior fixed point, we should further equate the two relevant forms of $\psi$ and $\varphi$.

Equating $\eqref{normal+interior fixed point+psi}$ and $\psi(z)=\frac{c_0^2-\alpha c_1}{c_0^2-\alpha c_1-(c_1-c_2)z}$, then comparing the constants and the coefficients of $z$ and $z^2$ respectively, we have that
\begin{align}
&c_0^2-\alpha c_1=\frac{(1-|p|^2\delta)(c_1-c_2)}{\overline{p}(1-\delta)},\label{6.3}\\
&\psi(p)(1-|p|^2)=1-|p|^2\delta,\psi(p)\neq0,c_0^2\neq\alpha c_1,c_1\neq c_2.\label{6.4}
\end{align}

(In fact, if $\psi(p)=0$ or $c_0^2-\alpha c_1$ or $c_1=c_2$, then $\psi$ is trivial.)

Further equating $\eqref{normal+interior fixed point+varphi}$ and $\varphi(z)=\frac{\alpha(\overline{\alpha}c_0^2-c_1)-(|\alpha|^2c_1-c_2)z}{\overline{\alpha}(c_0^2-\alpha c_1)-\overline{\alpha}(c_1-c_2)z}$, then comparing the constants and the coefficients of $z$ and $z^2$ respectively, we have that
$c_0^2-\alpha c_1=I_3(c_1-c_2),
\overline{\alpha}c_0^2-c_1=I_2I_3(c_1-c_2),
|\alpha|^2c_1-c_2=I_1(c_1-c_2).$
Combining \eqref{6.3} and \eqref{6.4}, we conclude that
\begin{align}
&c_0^2-\alpha c_1=I_3(c_1-c_2)\label{6.5}\\
&\overline{\alpha}c_0^2-c_1=I_2I_3(c_1-c_2)\label{6.6}\\
&|\alpha|^2c_1-c_2=I_1(c_1-c_2)\label{6.7}\\
&I_1=\frac{\overline{\alpha}(p^2-\delta)}{p(1-\delta)}, I_2=\frac{\overline{\alpha}p(1-\delta)}{\alpha(1-p^2\delta)},I_3=\frac{1-p^2\delta}{p(1-\delta)},\nonumber\\
&p\in\mathbb{R},c_0^2\neq\alpha c_1, \psi(p)\neq0, c_1\neq c_2, \delta\neq1, p\neq0\nonumber.
\end{align}
Substituting the expressions above into $\varphi(z)=\frac{\alpha(\overline{\alpha}c_0^2-c_1)-(|\alpha|^2c_1-c_2)z}{\overline{\alpha}(c_0^2-\alpha c_1)-\overline{\alpha}(c_1-c_2)z}$, we obtain $\eqref{normal+interior fixed point+varphi}$. This completes the proof.\hfill$\Box$

\begin{ex}\label{normal+C_2-symmetric+hyperbolic}
There is no $\mathcal{C}_2-$symmetric and normal weighted composition operator $W_{\psi,\varphi}$ on $H^2(\mathbb{D})$ if $\varphi$ is a hyperbolic linear-fractional self-map.
\end{ex}
\textit{Proof.} The result can be directly obtained Corollary \ref{normal+C_2-symmetric+aut} if $\varphi$ is a hyperbolic automorphism linear-fractional self-map and by the remark below Proposition 13 in \cite{BN2} if $\varphi$ is a hyperbolic non-automorphism linear-fractional self-map.

\begin{ex}\label{normal+C_2-symmetric+parabolic}
The weighted composition operator $W_{\psi,\varphi}:H^2(\mathbb{D})\rightarrow H^2(\mathbb{D})$ is $\mathcal{C}_2-$symmetric and normal when $\varphi$ is a parabolic linear-fractional self-map with the Denjoy-Wolff point $\zeta\in\partial\mathbb{D}$ if and only if
$$\psi(z)=\frac{c_0^2-\alpha c_1}{c_0^2-\alpha c_1-(c_1-c_2)z},$$
$$\varphi(z)=\frac{\overline{\alpha}\zeta^2(c_1-c_2)-(2\overline{\alpha}\zeta(c_1-c_2)-\overline{\alpha}(c_0^2-\alpha c_1))z}{\overline{\alpha}(c_0^2-\alpha c_1)-\overline{\alpha}(c_1-c_2)z},$$
$$(|\alpha|^2c_1-c_2+\overline{\alpha}(c_0^2-\alpha c_1))^2=4|\alpha|^2(c_1-c_2)(\overline{\alpha}c_0^2-c_1)$$ and Situation (i) or (ii) in Theorem \ref{normal+C_2-symmetric} is satisfied, where $c_0,c_1,c_2,\alpha\in\mathbb{D}\setminus\{0\},z\in\mathbb{D}.$
\end{ex}
\textit{Proof.} Since $\varphi$ is a hyperbolic automorphism linear-fractional self-map, then $\varphi$ has the Denjoy-Wolff point $\zeta=\frac{(|\alpha|^2c_1-c_2)+\overline{\alpha}(c_0^2-\alpha c_1)}{2\overline{\alpha}(c_1-c_2)}.$ Some trivial but tendious calaulstion shows that $\overline{\alpha}\zeta^2(c_1-c_2)=\alpha(\overline{\alpha}c_0^2-c_1).$ Note that
$$\zeta^2=\frac{(|\alpha|^2c_1-c_2+\overline{\alpha}(c_0^2-\alpha c_1))^2}{4\overline{\alpha}^2(c_1-c_2)^2}=\frac{\alpha(\overline{\alpha}c_0^2-c_1)}{\overline{\alpha}(c_1-c_2)}.$$
Hence, the result can be easily obtained by what we have observed above and some tedious calculation. This completes the proof.\hfill$\Box$

\end{document}